\documentclass[11pt,reqno]{amsart}
\usepackage{graphicx
}

\newcommand{\la}{\langle}
\newcommand{\ra}{\rangle}
\newcommand{\p}{\partial}

\newcommand{\Z}{\mathbb Z}

\newcommand{\R}{\mathbb R}

\newcommand{\im}{\operatorname{im}}

\renewcommand{\phi}{\varphi}
\newcommand{\arf}{\operatorname{Arf}}
\newcommand{\pd}{\operatorname{PD}}

\newcommand{\Hom}{\operatorname{Hom}}

\newtheorem{theorem}{Theorem}[section]
\newtheorem{lemma}[theorem]{Lemma}
\newtheorem{proposition}[theorem]{Proposition}

\theoremstyle{definition}
\newtheorem{definition}[theorem]{Definition}

\theoremstyle{remark}
\newtheorem{remark}[theorem]{Remark}

\numberwithin{equation}{section}

\begin{document}

\title[Knotted surfaces]{Topological triviality of smoothly knotted surfaces in $4$-manifolds}
\thanks{The second author was partially supported by NSF
Grant 0505605.}
\author[Hee Jung Kim]{Hee Jung Kim}
\address{Department of Mathematics \& Statistics\newline\indent
Hamilton Hall, Room 311\newline\indent McMaster
University\newline\indent 1280 Main Street West\newline\indent
Hamilton, Ontario L8S 4K1\newline\indent Canada}
\email{\rm{hjkim@math.mcmaster.ca}}
\author[Daniel Ruberman]{Daniel Ruberman}
\address{Department of Mathematics, MS 050\newline\indent Brandeis
University \newline\indent Waltham, MA 02454}
\email{\rm{ruberman@brandeis.edu}}

%\date\today

\begin{abstract}
Some generalizations of the Fintushel-Stern  rim surgery are known to
produce smoothly knotted surfaces.  We show that if the fundamental
groups of their complements are standard, then these surfaces are
topologically unknotted.
\end{abstract}
\maketitle

\section{Introduction}
Fintushel and Stern described~\cite{fs:knots} a surgery operation on
a torus $T$ embedded with trivial normal bundle in a smooth
$4$-manifold $X$, and used this construction to build many
$4$-manifolds with interesting properties.   The construction depends
on a choice of a knot $K$ in $S^3$, and the choice of a
diffeomorphism $\phi:  \p \nu(T)\to S^1\times \p E(K) $ (subject to
certain restrictions), and gives rise to a $4$-manifold $X_K(\phi)$.
When this operation is applied to a torus embedded in a neighborhood
of a surface $\Sigma$ in $X$ (and linking $\Sigma$ in a simple
fashion) the new manifold $X_K(\phi)$ is diffeomorphic to $X$.
However,  the construction, called rim surgery, produces a new
surface $\Sigma_K(\phi)$ in $X$ which Fintushel and Stern
showed~\cite{fs:surfaces} may be smoothly non-isotopic to $\Sigma$.

In the case that $\pi_1(X - \Sigma)$ is trivial, it is
straightforward to show that $\Sigma_K(\phi)$ is \emph{topologically}
isotopic to $\Sigma$ using topological surgery.  Recently,
S.~Finashin~\cite{finashin:surfaces} and the first-named
author~\cite{kim:surfaces} have used variations of rim surgery to
find smoothly knotted surfaces whose complements have (nontrivial)
cyclic fundamental groups.  It is interesting to ask whether these
surfaces are topologically nontrivial.  Indeed, Finashin~\cite[page
50]{finashin:surfaces} asks about the topological classification of
his surfaces, while the first author~\cite{kim:surfaces} showed that
for many knots $K$ and choices of $\phi$, the surface
$\Sigma_K(\phi)$ is topologically equivalent to the original
$\Sigma$.  In this paper, we show that any surface produced by a
torus knot surgery whose complement has cyclic fundamental group is
topologically standard.

To state the main result, we need to
define the surgery operation in question.
\begin{definition}\label{knotsurgery}
Let $Z$ be a 4-manifold containing a torus $T$ of self-intersection
$0$, and let $K$ be a knot with exterior $E(K)$.   Denote by $\mu_T$
the boundary of the normal disk of the torus, and let the
longitude/meridian of K be $\lambda_K$ and $\mu_K$ respectively.  Let
$\phi:  \p \nu(T)\to S^1\times \p E(K) $ be any
diffeomorphism such that $\phi_* \mu_T =\lambda_K $.\\
\indent(a)
\emph{Knot surgery}
is the operation $Z \to Z_K(\phi) = Z- (T \times D^2) \cup_\phi (S^1
\times E(K))$.\\
\indent(b) If $\pi_1  Z_K(\phi)\cong \pi_1
Z$ is cyclic, then the knot surgery will be called a \emph{cyclic
surgery}.
\end{definition}
Note that the boundary of $Z$ is not touched by this operation, so it
is meaningful to ask if there is a diffeomorphism or homeomorphism
$Z_K(\phi)\to Z$ extending the identity map on $\partial Z$.
Fintushel and Stern~\cite{fs:surfaces} investigated this construction
when $Z = X - \nu(\Sigma)$ is the exterior of an embedded surface in
the closed manifold $X^4$.  Gluing back in the neighborhood of the
surface gives a new embedding of $\Sigma$ in $X_K(\phi)$, with image
$\Sigma_K(\phi)$.  Fintushel and Stern focused particularly on a
torus $T$ that is the preimage in $\partial \nu(\Sigma)$ of a closed
curve $\alpha \subset \Sigma$, which they called a \emph{rim torus}.
In this case (referred to as \emph{rim surgery}) there is a canonical
identification of $X$ and $X_K(\phi)$, and we view $\Sigma_K(\phi)$
as lying in $X$.  If an appropriate Seiberg-Witten invariant of $X$
is non-trivial, and the Alexander polynomial of $K$ is not $1$, then
$\Sigma_K(\phi)$ is smoothly knotted.  On the other hand, if
$\pi_1(X-\Sigma) =1$, then rim surgery does not change the
topological type of the embedding.  Our main result is a
generalization of this last statement.

\begin{theorem}\label{T:unknotted}
Let $\Sigma \subset X$ be an embedded surface such that
$\pi_1(X-\Sigma) \cong \Z_d$.  Suppose that $T$ is a torus in the
complement of $\Sigma$, and that the knot surgery $(X,\Sigma)
\to(X_K(\phi),\Sigma_K(\phi))$ is cyclic. Then there is a pairwise
homeomorphism
$(X_K(\phi),\Sigma_K(\phi)) \to (X,\Sigma)$.
\end{theorem}
The proof of the theorem uses the traditional paradigm of surgery
theory: we will find a homotopy equivalence, modify it so that it is normally
bordant to the identity, and then show that the normal bordism may be
surgered to give an s-cobordism.

In certain cases, there is a canonical diffeomorphism between $X$ and
$X_K(\phi)$, and we regard $\Sigma$ and $\Sigma_K(\phi)$ as embedded
surfaces in $X$.   The examples we have in mind are the rim surgery
referred to previously, and the \emph{annulus rim surgery} in
Finashin's paper~\cite{finashin:surfaces}.   In this circumstance, it
is reasonable to ask if $\Sigma$ and $\Sigma_K(\phi)$ are isotopic.
\begin{theorem}\label{T:isotopy}
Suppose that $X$ is simply connected and $\pi_1(X- \Sigma)$ is
cyclic. If $\Sigma_K(\phi)$ is
obtained from $\Sigma$ by a cyclic rim surgery or annulus rim surgery, then
$\Sigma_K(\phi)$ is isotopic to $\Sigma$.
\end{theorem}
The final section of the paper contains some examples of cyclic
surgeries to which
Theorems~\ref{T:unknotted} and~\ref{T:isotopy} apply.
\begin{remark}
There are presumably alternate ways to obtain these results, most
notably the general methods for topological surgery with finite
fundamental group due to Hambleton and 
Kreck~\cite{hambleton-kreck:cancellation,hambleton-kreck:elliptic,hambleton-kreck:finite} that can perhaps
be applied in this context.
(Compare~\cite{finashin-kreck-viro:surfaces,kreck:surfaces}.)
We are able to use the more traditional surgery theory because we are
able to construct a specific homotopy equivalence between $Z$ and
$Z_K(\phi)$, whose normal invariant can be readily calculated by
geometric means.
\end{remark}
\section{Proofs of the theorems}\label{S:proofs}
We start with a simple way to recognize a homotopy equivalence
between $4$-manifolds with finite $\pi_1$.
\begin{lemma}\label{l:homotopic}
Let $M'$ and $M$ be oriented $4$-manifolds with finite fundamental
group. Suppose that $f: (M',\p M') \to (M,\p M)$ is a degree--$1$ map
with the following properties:
\begin{enumerate}
\item $f:\p M' \to \p M$ is a homotopy equivalence.
\item $f_*: \pi_1M'  \stackrel{\cong}{\longrightarrow} \pi_1M$.
\item $f_*: H_*(M';\Z)  \stackrel{\cong}{\longrightarrow} H_*(M;\Z)$.
\end{enumerate}

Then $f$ is a homotopy equivalence.
\end{lemma}
\begin{proof}
Denote by $\tilde M \to M$ the universal cover, and let $d$ be its degree.
By the Whitehead theorem (see~\cite{davis-kirk:at} for the statement
in the non-simply connected case) it suffices to show that the lift
$\tilde f$ of $f$ to the universal cover induces an isomorphism in
integral homology. Since the degree of $f$ is $1$, the same is true
for $\tilde f$, and so the induced map in homology is a surjection.
Now $H_1(\tilde M';\Z) = H_1(\tilde M;\Z) = 0$, and Poincar\'e
duality (plus the fact that $\tilde f$ is a homotopy equivalence on
the boundary) shows that $f_*:H_3(\tilde M';\Z)
\stackrel{\cong}{\longrightarrow} H_3(\tilde M;\Z)$.

Again, the fact that the degree of $\tilde f$ is one implies that
$$f_*:H_2(\tilde M';\Z)  {\longrightarrow} H_2(\tilde M;\Z),$$
a map between two free abelian groups, is a surjection.   On the other
hand, we know from (3) that the Euler characteristics of $M$ and $M'$
are equal, and so $\chi(\tilde M') = d\cdot \chi(M) =  d \cdot
\chi(M')  =\chi(\tilde M) $.  It follows that $H_2(\tilde M';\Z)$ and
$ H_2(\tilde M;\Z)$ have the same rank, and so $\tilde f_*$ must in
fact be an isomorphism.
\end{proof}
Let $X$ be a simply-connected oriented $4$-manifold.  For an
oriented embedded surface $\Sigma \subset X$, let $Z = X -  \nu(
\Sigma)$ be its exterior; note that $\pi_1(Z)$ is normally generated
by the meridian of $\Sigma$, which we will denote by $\mu_\Sigma$.
If $\pi_1(Z)$ is cyclic, then of course $\mu_\Sigma$ is a generator.
Suppose that $T \subset Z$ is an embedded torus such that $T\cdot T
= 0$, with a framing of its normal bundle, giving an identification
$\p\nu(T) \cong T^3$.  (In the case of rim surgery, the framing will
be canonical.)   Let $(S^3,K)$ be a knot with exterior $E(K)$.  The
choice of a diffeomorphism $\phi: T^3 \to S^1\times \p E(K)$ then
gives a manifold $Z_K(\phi) = Z -\nu(T) \cup_\phi (S^1\times E(K))$.

The gluing map $\phi$ can be encoded as a $3 \times 3$ matrix with respect to
a basis $\{\alpha,\beta,\mu_T\}$ for $H_1(T^3)$ and $\{S^1,
\mu_K,\lambda_K\}$ for $H_1(S^1 \times \p E(K))$.  We will always
assume that $\phi$ takes $\mu_{T}$ to $\lambda_K$, so the matrix will
have the form
\begin{equation}\label{matrix1}
\phi=\begin{pmatrix}
      a&b&0\\
      c&d&0\\
      p&q&1
\end{pmatrix}
\end{equation}
where $ad-bc=1$.

If $a=d=1$ and $b=q=0$ then such a gluing corresponds, in the
original setting of rim
surgery, to replacing the \emph{spinning} construction of
Fintushel-Stern with a combination of ($c$-fold) twist
spinning~\cite{zeeman:twist} and ($p$-fold) roll
spinning~\cite{fox:rolling,litherland:deform}.  The matrix
formulation of these classical constructions from high-dimensional knot theory is described in great detail in~\cite{plotnick:fibered}
and~\cite{montesinos:twins.I,montesinos:twins.II}.

\begin{lemma}\label{l:unknot}
Suppose that $\phi$ is a gluing map with $\phi(\mu_T)=\lambda_K$ ,
and that $K$ is the unknot $O$. Then $Z_O(\phi) \cong Z$ via a
diffeomorphism that is the identity on the boundary.
\end{lemma}
\begin{proof}
Consider the following  decompositions of $Z$ and $Z_O(\phi)$:
$$Z = Z- (T \times D^2) \cup_{id} (S^1 \times E(O))$$ and
$$Z_O(\phi) = Z- (T \times D^2) \cup_\phi (S^1 \times E(O)).$$ Then
it is sufficient to extend $id\circ\phi^{-1} :S^1 \times
\partial E(O) \to S^1 \times
\partial E(O)$ to a diffeomorphism $S^1 \times E(O) \to S^1 \times  E(O)$.
Using a diffeomorphism of $T^3$ that preserves the fiber $\mu_T$, we
may assume that~\eqref{matrix1} has the simpler form

\begin{equation}\label{matrix2}
\phi=\begin{pmatrix}
      1&0&0\\
      0&1&0\\
      p&q&1
\end{pmatrix}.
\end{equation}

Considering $S^1 \times E(O)$ as $S_{x}^1\times S_{y}^1\times D^2$,
we can write an element in  $S^1 \times E(O)$ as $(e^{2\pi
xi},e^{2\pi yi},re^{2\pi zi})$. Then we define an extension of
$\phi^{-1}$ as follows;
$$\phi^{-1}(e^{2\pi xi},e^{2\pi yi},re^{2\pi zi})=
(e^{2\pi xi},e^{2\pi yi},re^{2\pi (z-rpx-rqy)i})$$ where $0\le
x,y,z,$ and $r \le 1$.
\end{proof}

Because of Lemma~\ref{l:unknot}, it will suffice to show that for any
cyclic surgery involving the
$K$ there is a homeomorphism $Z_K(\phi) \cong Z_O(\phi)$.
\begin{lemma}\label{L:homotopy}
Let $K$ be a knot in $S^3$, and suppose that $Z \to Z_K(\phi)$ is a
surgery where $\pi_1(Z_K(\phi))$ is finite. There is an orientation preserving homotopy
equivalence $f: Z_K(\phi) \stackrel{\simeq}{\longrightarrow}
Z_O(\phi)$.
\end{lemma}
\begin{proof}
A straightforward obstruction theory argument produces a map $g:
E(K) \to E(O) \cong S^1 \times D^2$ that induces an isomorphism in
integral homology, takes meridian to meridian, and longitude to
longitude.  Take the product of this map with the identity map on
the circle, to get a homology equivalence $S^1 \times E(K) \to S^1
\times E(O)$.  Because we use the same gluing map $\phi$ for both
knots, this can be glued to the identity map on $Z- \nu(T)$ to get a
map $f: Z_K(\phi) \to Z_O(\phi)$.   The 5-lemma implies that $f$
induces an isomorphism on integral homology, it is an isomorphism on
$\pi_1$ as well.  The result follows by lemma~\ref{l:homotopic}
\end{proof}
We will refer to the homotopy equivalence provided by
Lemma~\ref{L:homotopy} as the \emph{canonical} homotopy equivalence.
It is automatically a normal map~\cite{wall:book}, and we would like
to construct a normal cobordism from $f$ to the identity map of
$Z_O(\phi)$, which will be denoted $Z'$ henceforth.  Recall that a
(topological) normal map $h: Y \to Z'
=Z_O(\phi)$, relative to the identity on the boundary, is classified by
an element $n(h) \in [Z',\partial Z'; G/TOP]$.  A standard
calculation~\cite{kirby-taylor:surgery} says that this is isomorphic to $H^4(Z',\partial Z')
\oplus H^2(Z',\partial Z';\Z_2)$, and that the first component of $n(h)$
is just $\frac18(\sigma(Y) - \sigma(Z'))$. In the case at hand, the
signature of $Y= Z_K(\phi) $ is the same as that of $Z'$, so we
concentrate on the second component, which we will denote by $S(f)$.

Note that the copy of $S^1 \times E(O)$ used in constructing $Z'$
contains a torus $T' = S^1 \times \mu_K$.  The homology class of this
torus in $Z'$ will be denoted $[T']$.
\begin{lemma}\label{split}
For the canonical homotopy equivalence $f$, the invariant $S(f)\in
H^2(Z',\partial Z';\Z_2)$ is given by $\arf(K)\pd(T')$.
\end{lemma}
\begin{proof}
Since  $H^2(Z',\partial Z';\Z_2) \cong \Hom (H_2(Z',\partial
Z';\Z_2); \Z_2)$, the class $S(f)$ is determined by its evaluation on
(possibly nonorientable) surfaces $(G,\partial G) \subset
(Z',\partial Z')$.  The following recipe for this evaluation is generally viewed as part of Sullivan's characteristic variety theorem~\cite{sullivan:gt,sullivan:hauptvermutung}.  We could not find this exact statement in the literature, and so have supplied a proof in an appendix to this paper.
Make $f$ transverse to $G$, to get a 2-dimensional
normal map $f^{-1}G \to G$.  Then $\langle S(f), [G,\partial
G]\rangle $ is the surgery obstruction of this normal map (in
$L_2({1}) \cong \Z_2$) which in turn is given by the Arf invariant.
Note that this evaluation depends on the relative homology class
$[G,\partial G]$, but not on the choice of surface representing that class.

Now given a surface $G \subset Z'$, its mod-$2$ intersection number
with $T'$ is either $0$ or $1$.  By changing $G$ by a homology (in
fact by an embedded cobordism) we  may assume that the intersection
of G with $S^1 \times E(O) \subset Z'$ is, correspondingly, either
empty or a single copy of the disk spanning $O$ in $E(O)$.  This
uses the fact that the boundary of this disk is preserved by the
gluing map $\phi$.  If $G\cdot T' = 0$, so the intersection is empty,
then $f$ is a homeomorphism on $f^{-1}G$, and so the Arf invariant
computing $\langle S(f), [G,\partial G]\rangle $ must also vanish.

On the other hand, if  $G\cdot T' = 1$, so the intersection is the
meridional disk, then we can decompose $G = G_0 \cup D^2$ where the
$D^2$ is the disk bounding $O$ in $E(O)$.  Clearly $f^{-1}G$
decomposes in a corresponding way, where $f: f^{-1}G_0 \to G_0$ is a
homeomorphism, and $f^{-1}D^2$ is a a Seifert surface for $K$,
mapping with degree $1$ onto $D^2$.   Thus $\langle S(f), [G,\partial
G]\rangle $ is the Arf invariant of a quadratic form defined on $\ker
[H_1( f^{-1}G;\Z_2) \to H_1(G;\Z_2)]$.  But this kernel is clearly the
same as the (mod 2) homology of the Seifert surface, and it is
straightforward to identify the quadratic form with the one that
gives the Arf invariant of $K$.  Putting together these two cases, we
get that $\langle S(f), [G,\partial G]\rangle = \arf(K) G \cdot T' =
\arf(K)\langle \pd(T'),[G,\partial G]\rangle $.
\end{proof}
     From this lemma, if the Arf invariant of $K$ is trivial, then the
homotopy equivalence constructed in Lemma~\ref{L:homotopy} is
normally cobordant to the identity.   On the other hand, if $\arf(K)
=1$, then $g$ is not normally cobordant to the identity, and so this
construction does not work.   However, we will see that this Arf
invariant can be `absorbed into the complement' by choosing a
different homotopy equivalence.
\begin{proof}[Proof of Theorem~\ref{T:unknotted}]
We claim that there is a homotopy equivalence (relative to the
boundary) $f' : Z' \to Z'$ with the same normal invariant as $f$.
If the Arf invariant of $K$ is trivial, then the normal invariant of
$f$ is trivial, and so $f$ is normally cobordant to the identity.  In
fact, the normal cobordism is easy to construct `by hand'.
Following~\cite{freedman:concordance,freedman-quinn}, the vanishing
of $\arf(K)$ gives rise to an explicit normal cobordism, relative to
the boundary, of the map $g$ to the identity of $E(O)$.   As in the
proof of Lemma~\ref{L:homotopy}, this normal cobordism (crossed with
$S^1$) can be glued to the identity map on $S^1 \times (Z- \nu(T))$
to give the desired normal cobordism.

If $\arf(K) =1$, then we make use of a classic construction of
Wall~\cite[\S16]{wall:book} that is discussed carefully
in~\cite[section 5]{cochran-habegger:homotopy}.   Given a class in
$\pi_2(Z')$ represented by a map $\alpha: S^2 \to Z'$, consider the
composition $\alpha\circ \eta^2\in \pi_4Z'$, where $\eta^2$ is the
generator of $\pi_4(S^2)$.  Then the following composition, say
$f_\alpha$, is a homotopy equivalence (rel $\partial Z'$):
$$
Z'  \longrightarrow  Z' \vee S^4  \stackrel{\alpha\circ
\eta^2}{\longrightarrow} Z'.
$$
The normal invariant $S(f_\alpha)$ is computed geometrically in
Theorem 5.1 of~\cite{cochran-habegger:homotopy}, and is given by
\begin{equation}\label{ch}
(1 + \langle w_2(Z'),\alpha\rangle)\pd(\alpha).
\end{equation}
Here we have used the same notation for the homotopy class $\alpha$
and the homology class it carries.

Note that the image of the Hurewicz map is determined via the Hopf
exact sequence
$$
\pi_2(X) \stackrel{h}{\longrightarrow} H_2(X){\longrightarrow}
H_2(\pi_1(X)) \to 0.
$$
Since $H_2(\Z_d) = 0$, the Hurewicz map is onto, so in particular the
class $T'$ is spherical and we get a homotopy equivalence $f_{T'}: Z'
\to Z'$.  Recall that $T' = S^1 \times \mu_K \subset S^1\times E(K)$,
so it has trivial normal bundle.  Thus $\langle w_2(Z'),T'\rangle
\equiv T' \cdot T' \equiv 0\pmod{2} $ and so by Lemma~\ref{split} and
formula~\ref{ch},  both $S(f_{T'})$ and $S(f)$ are given by
$\pd(T')$.  This establishes the claim.

A normal cobordism $W$ from $f:Z_K(\phi)\to Z'$ to the homotopy
equivalence $f':Z' \to Z'$ has a surgery obstruction that \emph{a
priori} lies in the Wall group $L_5^h(\Z[\Z_d])$. But in
fact~\cite[section 11]{hambleton-taylor:guide} this group, and the
group $L^s(\Z[\Z_d])$ vanish, so that $W$ may be surgered to give an
s-cobordism rel boundary between $Z_K(\phi)$ and $Z'$, which is a
topologically a product.  Hence the identity map between $\partial
Z_K(\phi)$ and $\partial Z'$ extends to a homeomorphism.   (For $d$
odd and $Z$ closed, this argument can be found
in~\cite{lawson:inertial}.)   Finally, since the map on the boundary
is the identity, the homeomorphism extends over $X$.
\end{proof}
\begin{proof}[Proof of Theorem~\ref{T:isotopy}]
Note that when $T$ is a rim torus (or annular rim torus as
in~\cite{finashin:surfaces}),
the homology class $T$ is trivial
in $H_2(X)$.  So the homeomorphism constructed in the proof of
theorem~\ref{T:unknotted} has trivial normal invariant when viewed as
a map from $X$ to $X$.  Hence by~\cite{quinn:isotopy} (again
see~\cite{cochran-habegger:homotopy} for full details) that
homeomorphism is homotopic to the identity.  Then the work of
Perron~\cite{perron:isotopy2} and Quinn~\cite{quinn:isotopy} says
that this homeomorphism is isotopic to the identity. The isotopy
takes $\Sigma$ onto $\Sigma_K(\phi)$.
\end{proof}
\section{Examples of cyclic surgeries}
\label{S:examples}
In this section, we discuss those examples we know of gluing maps
which produce a surface whose complement has cyclic fundamental
group. In many cases, these will be smoothly knotted, but
Theorem~\ref{T:isotopy} shows that they are topologically standard.
In all of these, $\Sigma$ is a surface in the simply-connected
manifold $X$ with $\pi_1(X-\Sigma)\cong \Z/d$.

\subsection{Rim Surgery~\cite{fs:surfaces}}
   Let $T$ be a rim
torus in $X-\Sigma$, i.e. of the form
$\alpha\times\mu_\Sigma$ where $\alpha$ is a non-separating curve in
$\Sigma$.  The simplest gluing map  $\phi: \p\nu(T) (\cong T\times\p
D^2) \to S^1\times \p E(K) $ sends $\alpha\mapsto [S^1]$,
$\mu_{\Sigma}\mapsto \mu_K$, and $\p D^2\mapsto \lambda_K$.  As shown
in~\cite[Example 3.2]{kim:surfaces}, the fundamental group of $X -
\Sigma_K(\phi)$ contains the fundamental group of the $d$-fold
branched cover of $(S^3,K)$, and hence (by a strong form of the Smith conjecture~\cite{morgan-bass:smith}) is bigger than $\Z/d$.

Cyclic surgeries arise when we perform twisting ($m$ times) and
rolling ($n$ times) in the gluing map $\phi$, which has the form:
\begin{equation}\label{matrix3}
   \phi=\begin{pmatrix}
      1&0&0\\
      m&1&0\\
      n&0&1
\end{pmatrix}
\end{equation}

%%%%%%%
The terminology (twisting and rolling) is explained in Section~\ref{S:proofs}.
In the case of $m$-twisted rim surgery (i.e. $n =0$) it is shown in
~\cite{kim:surfaces} that for any knot $K$,
$\pi_1(Z_K(\phi))$ is cyclic when $d\equiv\pm 1 (\mbox{mod\ } m)$,
but can also be non-cyclic when this condition is violated.   For instance, if $d | m$ and $K$ is non-trivial, then as in Example 3.2 of~\cite{kim:surfaces}, the group $\pi_1(Z_K(\phi))$ will not be cyclic.   We  will show that $\phi$ is a
cyclic surgery for any knot $K$ in a more general case which includes an arbitrary amount of rolling.

\begin{proposition}\label{fundamentalgroup} If $(d,m)= 1$, and $n$ is arbitrary,
then in the above language, the $m$-twisted $n$-rolled rim surgery has
$\pi_1(Z_K(\phi)) = \Z_d$ for any knot $K$.
\end{proposition}

To prove this, we will first describe $m$-twists and $n$-rolls of
rim surgery following~\cite{kim:surfaces}. Given a non-separating
curve $\alpha$ in $\Sigma$, choose a neighborhood of $\alpha$ in $X$
of the form $S^1\times I\times D^2=S^1\times B^3$, where
$S^1\times I$ is a neighborhood of $\alpha$ in $\Sigma$. Adjusting a trivialization of the normal bundle $\nu(\Sigma\times I )$ in $X$, we may assume that the push off of the curve $\alpha$ along the trivialization is homologically trivial in $X-\Sigma$.  Now we define
self diffeomorphisms denoted by $\tau$ and $\rho$ on $(S^3,K)$ as it
follows. Define a twist, $\tau$, by
\begin{equation}\label{tau}
\tau(\overline\theta, e^{i\varphi}, t) =(\overline\theta, e^{i(\varphi + 2\pi
t)}, t) \quad \mbox{for} \quad
(\overline\theta, e^{i\varphi}, t)\in K\times{\partial{D^2}}\times{I}
\end{equation}
and otherwise, $\tau(y)=y$. (Here we use $K\cong S^1\cong \R/\Z$.)\\
Similarly, a roll, $\rho$, is given by
\begin{equation}\label{rho}
\rho(\overline\theta, e^{i\varphi}, t)  =(\overline{\theta+t}, e^{i\varphi}, t) \quad \mbox{for} \quad
(\overline\theta, e^{i\varphi}, t)\in K\times{\partial{D^2}}\times{I}
\end{equation}
and otherwise, $\rho(y)=y$.  It is useful to note that $\rho$ and $\tau$ commute.

Divide $(S^3,K)$ into two arcs $ (B^3,K_+) \cup (B^3_-,K_-)$ where the second pair is unknotted, and $B^3_-$ lies inside the tubular neighborhood $K \times D^2$.  The pair $(X,\Sigma_K(\phi))$ defined  by $m$-twisted and $n$-rolled
rim surgery is then obtained by taking out the neighborhood of
$\alpha$ and gluing back the mapping torus of $(B^3,K_+)$ with monodromy given by powers of the
diffeomorphisms $\tau$ and $\rho$.
\begin{equation}\label{twist-rolling}
(X,\Sigma_K(\phi))=(X,\Sigma) - S^1\times(B^3,I)\cup_{\partial}
S^1\times_{\rho^n\tau^m} (B^3,K_{+})
\end{equation}

As in~\cite{kim:surfaces}, we observe that this description of $(X,\Sigma_K(\phi))$ is
equivalent to the one obtained by performing knot surgery in
Definition~\ref{knotsurgery} along the rim torus $T$ and the given matrix $\phi$ in
~\eqref{matrix3}.

Now let's consider $d$-fold covers of $X$ branched along $\Sigma$
and $\Sigma_K(\phi)$ denoted by $Y$ and $Y_K(\phi)$ respectively.
From the decomposition of $(X,\Sigma_K(\phi))$ in
~\eqref{twist-rolling} and the choice of the curve $\alpha$, we easily observe that $Y_K(\phi)$ is obtained by doing a surgery on $Y$ as follows:
\begin{equation}\label{branchcover1}
(Y_K(\phi),\Sigma_K(\phi)) =(Y,\Sigma)-S^1\times
(B^3,I)\cup_{\partial} S^1 \times_{\tilde\rho^n\tilde\tau^m}
(B^3,K_{+})^d
\end{equation}
where $(B^3,K_{+})^d$ is a $d$-fold cover of $B^3$ branched along
$K_{+}$ and $\tilde\rho$,$\tilde\tau$ are lifts of $\rho$,$\tau$
respectively. Here we consider $(B^3,K_{+})^d$ as the punctured
$d$-fold branched cover $(S^3,K)^d$ of $(S^3,K)$. The lift
$\tilde\rho$ of $\rho$ into $(S^3,K)^d$ is a map defined by a
rolling along the lifted knot of $K$ on a collar of
$\partial{E(K)}^d$ in a $d$-fold cover of the exterior ${E(K)}^d$ as
described in ~\eqref{rho} and otherwise, the identity. However
$\tilde\tau$ is a little more complicated.  Let $\sigma$ be the
canonical covering translation of the group $\Z_d$ of covering
transformations that is a rotation by $2\pi/d$ about the branch set.
Then the following map gives the lift of $\tilde\tau$ on
$(S^3,K)^d$.

\begin{equation}\label{lifttau}
\tilde\tau =
\begin{cases}
\sigma(x) & \text{if $x\in {E(K)}^d-\partial {E(K)}^d\times I$}\\
(\overline{\theta}, e^{i(s/d\cdot 2\pi +\varphi)},s) & \text{if
$x=(\overline{\theta}, e^{i{\varphi}},s)\in \partial {E(K)}^d\times
I$}\\
 x& \text{otherwise}
\end{cases}
\end{equation}

Note that $\tilde\rho$ is isotopic to the identity on
$(S^3,K)^d$. The isotopy $F_t$ between identity and
$\tilde\rho^{-1}$ is defined on a collar of $\partial{E(K)}^d$ as
follows:
$$F_t(\overline{\theta}, e^{i{\varphi}},s)
=(\overline{\theta-ts}, e^{i\varphi}, s) \quad \mbox{for} \quad
(\overline{\theta}, e^{i{\varphi}},s)\in
\partial{E(K)}^d\times I$$
This isotopy induces a homeomorphism
$(S^3,K)^d\times_{\tilde\rho}S^1\to (S^3,K)^d\times S^1$. So we
rewrite $Y_K(\phi)$ with a gluing map $f$ along the boundary.
\begin{equation}\label{branchcover}
(Y,\Sigma)-(B^3,I)\times
S^1\cup_{f^n}(B^3,K_+)^d\times_{{\tilde\tau}^m}S^1
\end{equation}

Here we need to make use of some results from Plotnick's paper ~\cite{plotnick:fibered}.  Consider a plumbing $P$, at two points, of two copies of $S^2\times D^2$.  The pair of cores of the $S^2\times D^2$s is called a `twin' by Montesinos~\cite{montesinos:twins.I,montesinos:twins.II}; note that $\partial P$ is a $3$-torus.  Plotnick~\cite{plotnick:fibered} constructed infinitely
many homotopy spheres by gluing to $ S^1\times E(K)$ using an identification $A$
defined on the boundary of $P$. More explicitly, write
$$\Omega_A=P\cup_{A} S^1\times E(K)$$
where $A$ is expressed by a matrix form according to a certain basis
$\{e_1,e_2,e_3\}$ on $H_1 (\partial P)$ and each $e_i$ represents a curve on $\partial P$.
(For details , see~\cite{plotnick:fibered}.) The general
form of $A$ for which $\Omega_A$ is a homology sphere  is
\begin{equation}\label{matrix4}
   \begin{pmatrix}
      p&k&0\\
      -\gamma&\beta&0\\
      -\alpha\gamma+bp&\alpha\beta+bk&1
\end{pmatrix}.
\end{equation}

In many cases, this construction produces $S^4$. In $\Omega_A$, we
get a knot  $A(K)$ determined by one of the two cores of $P$; in the description in ~\cite{plotnick:fibered} it is the second of the spheres.   An interesting result is that the knot $A(K)$ in $\Omega_A$
 is fibered. Moreover, in some cases the fiber is
described explicitly in terms of branched covers ~\cite[Theorem 5.6]{plotnick:fibered}.

\begin{proof}[Proof of Proposition~\ref{fundamentalgroup}]

Consider $(B^3,K_+)^d\times_{{\tilde\tau}^m}S^1$ in the decomposition ~\eqref{branchcover}  of  $Y_K(\phi)$.  Then there is an interesting connection with Plotnick's construction ~\cite{plotnick:fibered}. To see this, consider
our assumption $(d,m)= 1$. So there are $\gamma$ and $\beta$ such
that $d\gamma+m\beta=1$.  If we do Plotnick's construction using the following matrix
\begin{equation}\label{matrix5}
   A=\begin{pmatrix}
      m&d&0\\
      -\gamma&\beta&0\\
      0&0&1
\end{pmatrix}
\end{equation}
then the resulting space $\Omega_A$ is the $m$-fold cyclic branched cover of the $d$-twist spin of $K$
~\cite{zeeman:twist}. Theorem 5.6 in ~\cite{plotnick:fibered} says
that the knot $A(K)$ produced in $\Omega_A$ is fibered and its fiber
is exactly $(B^3,K_+)^d$. Moreover its characteristic map is same as
${\tilde\tau}^m$ described in ~\eqref{lifttau}. In other words,
$(B^3,K_+)^d\times_{{\tilde\tau}^m}S^1$ is the complement of  $A(K)$
in $\Omega_A$.

Comparing this gluing map A ~\eqref{matrix5}  with the general form ~\eqref{matrix4},
in our circumstance we have $p=m$, $k=d$, $-\alpha\gamma+bm=0$ and
$\alpha\beta+bd=0$. Since $d\gamma+m\beta=1$, we have
$d\gamma{b}+m\beta{b}=b$. Since $bm=\alpha\gamma$ and
$bd=-\alpha\beta$, $-\alpha\beta\gamma+\alpha\beta\gamma=b$. So
$b=0$. This implies that $\alpha=0$.
 According to Corollary 6.1 in
~\cite{plotnick:fibered}, this means that $\Omega_A$ is
smoothly $S^4$. Thus, $(B^3,K_+)^d\times_{{\tilde\tau}^m}S^1$ is
smoothly $S^4-A(K)$.   Then $Y_K(\phi)$ is isomorphic to
$$(Y,\Sigma)-(B^3,I)\times S^1\cup_{f^n} S^4-A(K).$$
This implies that $Y$ and $Y_K(\phi)$ are homeomorphic.

Now we consider $d$-fold unbranched covers of $Z=X-\Sigma$ and
$Z_K(\phi)=X-\Sigma_K(\phi)$ denoted by $\widetilde{Z}$ and
$\widetilde{Z_K(\phi)}$ respectively. Considering the description ~\eqref{branchcover1} of the branched cover
$Y_K(\phi)$, we can similarly write the
unbranched cover $\widetilde{Z_K(\phi)}$ as
$$\widetilde{Z}-S^1\times (B^3,I)\cup_{\partial}
S^1 \times_{\tilde\rho^n\tilde\tau^m} E(K)^d.$$

Let's apply Van Kampen Theorem to this decomposition. Since  $\widetilde{Z}-S^1\times (B^3,I)$ is homotopy equivalent to $\widetilde{Z}$ and $\pi_1\widetilde{Z}=\{1\}$,  the fundamental
group $\pi_1(\widetilde{Z_K(\phi)})$ is isomorphic to
$\pi_1(S^1 \times_{\tilde\rho^n\tilde\tau^m} E(K)^d)/\im j$ where $j:\pi_1( S^1\times (S^2-\{2 pts\}))\to \pi_1(S^1 \times_{\tilde\rho^n\tilde\tau^m} E(K)^d)$ is the inclusion  homomorphism.
Eventually, we have
\begin{equation}\label{pi_1:unbranchcover}
\pi_1(\widetilde{Z_K(\phi)})\cong\la \pi_1(E(K)^d,*) \mid \mu_{\tilde K}=1,
\beta =\tilde\rho_{*}^{n}\tilde\tau_{*}^{m}(\beta),
\forall\beta\in\pi_1(E(K)^d,*)\ra
\end{equation}
 where $\mu_{\tilde K}$ is a
meridian of the lifted knot $K$. To assert $\pi_1(Z_K(\phi)) = \Z_d$, it is sufficient to show that
$\pi_1(\widetilde{Z_K(\phi)})$ is trivial.
First, note that $\pi_1(Y_K(\phi))$ is trivial.
Applying Van Kampen Theorem to the decomposition for $Y_K(\phi)$ in ~\eqref{branchcover1} , we rewrite the group $\pi_1(Y_K(\phi))$ as
$$\la \pi_1((B^3,K_+)^d,*) \mid
\beta =\tilde\rho_{*}^{n}\tilde\tau_{*}^{m}(\beta),
\forall\beta\in\pi_1((B^3,K_+)^d,*)\ra .$$

Since $\pi_1((B^3,K_+)^d,*)\cong \pi_1(E(K)^d,*)/\la \mu_{\tilde
K}\ra$, this implies that the generator and relations in ~\eqref{pi_1:unbranchcover} makes
$\pi_1(\widetilde{Z_K(\phi)})$ trivial as well.

\end{proof}
\begin{remark}   We were led to the statement of Proposition~\ref{fundamentalgroup} by some computations with the group theory program GAP~\cite{GAP4}.  We
investigated $m$-twisted, $n$-rolled rim surgery for some simple
knots.    We found, for instance, that if $K$ is a figure 8 knot,
then $1$-twisted, $n$-rolled rim surgery is a cyclic surgery for
$n\leq 100$ (independent of $d$).   However, the computations did not
finish in a reasonable amount of time for more complicated knots and
larger values of $m$ or $n$.  In retrospect, it seems remarkable that Plotnick was able to show that in many cases $\Omega_A$ is actually $S^4$, whereas showing that its fundamental group is trivial seems somewhat difficult, computationally.
\end{remark}
\subsection{Finashin's Construction}\label{finashin}
Let's review briefly the annulus rim surgery construction
from~\cite{finashin:surfaces}. Suppose that there is a
smoothly embedded surface $M$ in $X$, called a `membrane', such that
$M\cong S^1\times I$, $M\cap \Sigma=\partial M$ and $M$ meets $\Sigma$ normally along $\partial M$. By adjusting a trivialization
of its regular neighborhood $U$, we can assume that $U(\cong
S^1\times B^3)\cap \Sigma=S^1\times f$, where $f=I_0\sqcup
I_1=I\times\partial I$ is a disjoint union of two unknotted segments
of a part of the boundary of a band $b=I\times I$ in $B^3$. Here the
band $b=I\times I$ is trivially embedded in $B^3$ and the
intersection $I\times I\cap
\partial B^3=\partial I\times I$.
Then given a knot $K$ in $S^3$, we consider a band $b_K\subset{B^3}$
obtained by knotting $b$ along $K$ and denote by $f_K$ the pair of
arcs bounding $b_K$. Note that the framing of $b_K$ is chosen the
same as the framing of $b$ (See ~\cite[Fig. 1]{finashin:surfaces}).
Then we get a new surface $\Sigma_{K}$ as follows:
$$(X,\Sigma_K)=(X,\Sigma) -S^1\times (B^3,f)\cup S^1\times (B^3,f_K)$$

There is an equivalent description for this construction using knot
surgery. Let
$\gamma$ be a meridian of $b$ in $B^3$. Then there is a torus
$T\subset \nu(M)$ corresponding to $S^1\times\gamma\subset S^1\times
(B^3,f)$, and $\Sigma_K$ is constructed by knot surgery on this
torus.  In Finashin's
original construction, the gluing map $\phi$ sends  $[S^1]\mapsto
[S^1]$, $\gamma\mapsto
\mu_K$, and $\p D^2\mapsto \lambda_K$.   More generally, we can do
$m$-twists and
$n$-rolls in this construction, resulting in a gluing map $\phi$  of  ~\eqref{matrix3}.
\begin{proposition}
Knot surgery using $\phi$ of ~\eqref{matrix3} preserves the fundamental group
$\pi_1(Z_K(\phi)) = \Z_d$. In other words, $m$-twisted and $n$-rolled
annulus rim surgery is a cyclic surgery.
\end{proposition}
\begin{proof}
Let's consider the process of
$m$-times twist-spinning and $n$-times rolling of a knot $K$ in
$S^3$. Choose a regular neighborhood $\nu (K)$ containing a knotted
band $b_K$ described in Finashin's construction.  We define self diffeomorphisms $\tau$
and $\rho$ on $(S^3, \nu (K))$ in the same way as~\eqref{tau} and ~\eqref{rho}.
Then $\tau$ is fixed on $\p E(K)$ and with support in a
collar of $\p E(K)$ and  induces a conjugation by a meridian.
Similarly, $\rho$ is a conjugation by a longitude. Then we can
observe that  $Z_K(\phi)$ obtained by knot surgery (Definition~\ref{knotsurgery}) along the above torus $T$  is equivalent to
$$Z-S^1\times (B^3,f)\cup S^1\times_{\rho^n\tau^m} (B^3-f_K).$$
Applying the Van Kampen
theorem to this decomposition of $Z_K(\phi)$
and following the similar argument in
~\cite{finashin:surfaces}, we can get that
$$\pi_1(Z_K(\phi))\cong \pi_1(S^1\times_{\rho^n\tau^m}
(B^3-f_K))/j(k)$$ where $j:\pi_1(S^1\times \p (B^3-f)\cong S^1\times
(S^2-\{4pts\}))\to \pi_1(S^1\times_{\rho^n\tau^m} (B^3-f_K))$ is the
inclusion homomorphism and $k$ is the kernel of the inclusion
homomorphism $\pi_1(S^1\times \p (B^3-f))\to \pi_1(Z-S^1\times
(B^3,f))$. The group $\pi_1(S^1\times \p (B^3-f)\cong S^1\times
(S^2-\{4pts\}))$ is generated by $[S^1]$ and three generators
denoted by $a_1$, $a_2$, and $a_3$ for $S^2-\{4pts\}$ where
$i(a_1)$= $i(a_2)$ is represented by a loop around one arc of $f$ in
$B^3$ and $i(a_3)$ is represented by a loop around $f$ under the
inclusion $i:\pi_1( \p (B^3-f))\to \pi_1(B^3-f)$ (See Fig. 1).
\begin{figure}[!ht]
\begin{center}
\includegraphics[scale=.8]{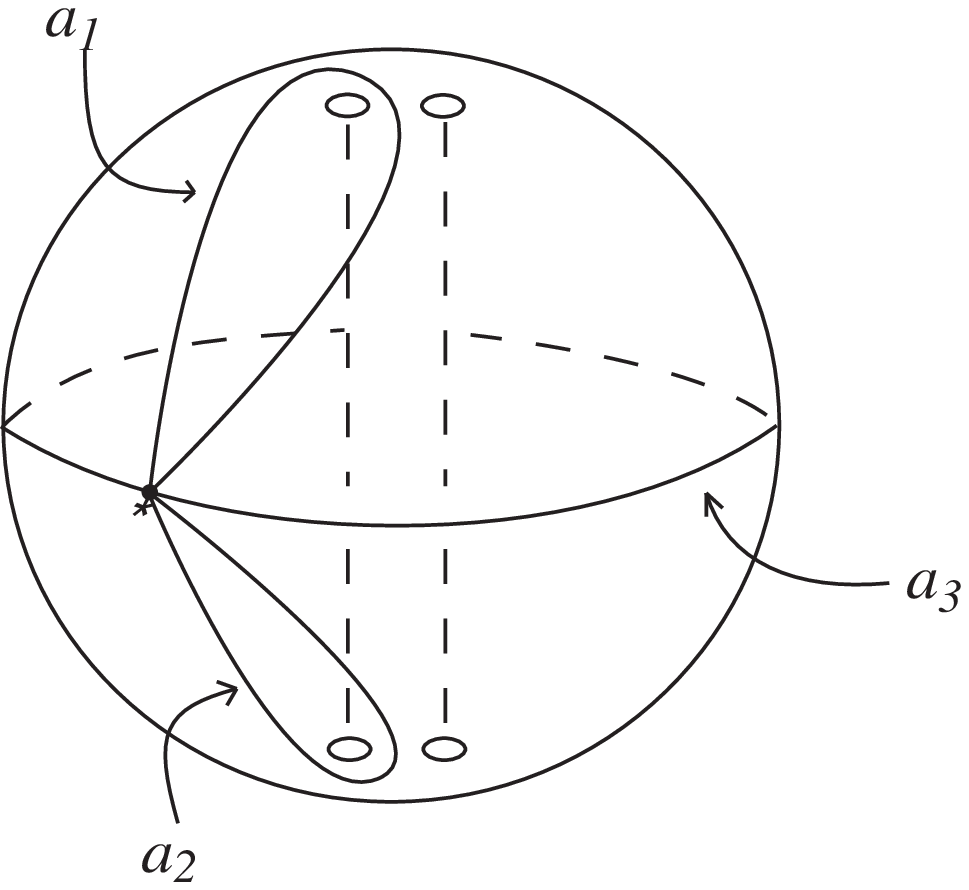}
\caption{}
\end{center}
\end{figure}
So the kernel $k$ is generated by $[S^1]$, $a_1^d$, $a_1a_2^{\pm 1}$,
and $a_3$. Since $j([S^1])$ kills the first factor of
$\pi_1(S^1\times_{\rho^n\tau^m} (B^3-f_K))$, the group
$\pi_1(Z_K(\phi))$ can be written by
$$\la \pi_1(B^3-f_K,*) \mid a_1^d=1,
\beta =\rho_{*}^{n}\tau_{*}^{m}(\beta), a_3=1, a_1a_2^{\pm 1}=1,
\forall\beta\in\pi_1(B^3-f_K,*)\ra $$ where we omit $j$ to simplify
the notation.  As discussed in ~\cite{finashin:surfaces}, the
relation $a_3$ has the same effect on $\pi_1(B^3-f_K,*)$ as
attaching a $2$-cell along a loop $m_b$ turning around $b_K$ in
$B^3$. Adding this $2$-cell makes $f_K$ an unknotted arc $I$ in
$B^3$ (See Fig. 2 which is the same figure  in
~\cite{finashin:surfaces}).
\begin{figure}[!ht]
\begin{center}
\includegraphics[scale=.8]{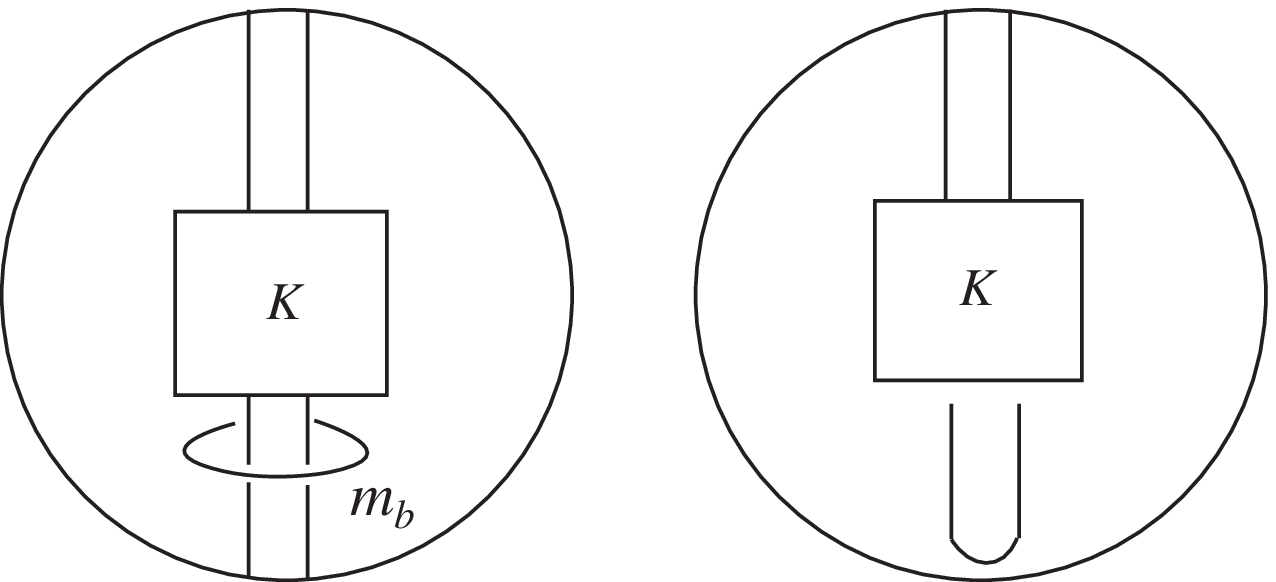}
\caption{}
\end{center}
\end{figure}
So only one generator $a_1$ of $\pi_1( \p (B^3-f))$ becomes the
generator for $\pi_1(B^3-I)$ (to see this, note that $a_1$ is
represented by a loop around $I$). By geometric observation, we know
that $\tau_{*}(a_1)$ is the conjugate of $a_1$ by $a_3$ and
$\rho_{*}(a_1)$ is the conjugate of $a_1$ by a longitude
$\lambda_K$. But since $a_3$ and the longitude are killed in the
group $\pi_1(B^3-I)$, the factorization by $j(k)$ makes the group
$\pi_1(Z_K(\phi))$ into a cyclic group $\Z_d$.
\end{proof}

\appendix
\section{Normal maps in dimension $4$}
In the proof of Lemma~\ref{split}, we gave a formula for evaluating the normal invariant of a map of $4$-manifolds.   It was pointed out to us by Ian Hambleton that Sullivan's characteristic variety theorem, on which this evaluation is based, does not have a complete proof in the literature.   Because $4$ is a fairly low dimension, it is possible to assemble a straightforward proof based on well-established facts in surgery theory, and we present such a proof in this appendix.

We treat the closed PL case first; extension to the bounded case is routine.  The TOP case follows by the same reasoning using the computation~\cite{kirby-siebenmann} of the homotopy groups of $G/TOP$. The statement we want is that PL normal maps to a closed oriented manifold $X$ are determined by two invariants, $S(f)$ and $\sigma(f)$ where $\sigma$ is the simply-connected surgery obstruction, given by the difference in signature divided by $8$.   By the Universal Coefficient Theorem, $S(f) \in H^2(X;Z/2)$ is determined by the pairings $\la S(f),S\ra$  for $S$ a surface (not necessarily orientable) in X.   The claim is that (for $f$ transverse to $S$) this evaluation is the  Arf invariant of the normal map $f^{-1}(S) \to S$.

The proof starts from the fact, proved by transversality and valid in any dimension, that PL normal maps to a closed manifold X are given by homotopy classes $[X,G/PL]$.  Now it is well-known (cf. for example~\cite[\S4]{madsen-milgram:surgery}) that the first two nonzero homotopy groups of $G/PL$ are $\pi_2=\Z/2$ and $\pi_4=\Z$, and $\pi_5=0$.  The first $k$-invariant in the PL case is $\delta Sq^2$; in the TOP case it is trivial.  This is enough to determine that $[X,G/PL]$ is the set $\{(a,s) \in H^2(X;\Z/2) \oplus H^4(X;\Z) \,|\, a^2=s \pmod 2 \}$.

Suppose that $F: M \to X$ is a normal map, and that $S \subset X$ is an embedded surface with $F$ transverse to $S$.  Then the restriction of $F$ is a normal map $f: F^{-1}(S) \to S$.  Moreover, the classifying map for $f$ = normal map of surfaces in $[S,G/PL]$ is the composition of the classifying map for $F$ in $[X,G/PL]$ with the inclusion $i:S \to X$.  This holds because the normal bundle of 
$F^{-1}(S)$ in $M$ is the pull-back of the normal bundle of $S$ in $X$.

The formula we need then follows from the following statement: if $f: S' \to S$ is a normal map of surfaces classified by a map $g: S \to G/PL$, and $a\in H^2(G/PL;\Z/2) = \Z/2$ is the generator, then $g^*(a) = \arf(f)$.   This can be seen in many ways, e.g. it follows directly from~\cite[Theorem 4.1]{rourke-sullivan:arf}. 

\section*{Acknowledgments}
We thank Ian Hambleton for some helpful correspondence on surgery theory, and for urging us to include a proof of our assertions about the normal invariant.
%\bibliography{mathbib}

\begin{thebibliography}{10}

\bibitem{cochran-habegger:homotopy}
Tim~D. Cochran and Nathan Habegger, \emph{On the homotopy theory of simply
  connected four manifolds}, Topology \textbf{29} (1990), no.~4, 419--440.
  \MR{\bf 1071367 (91h:57006)}

\bibitem{davis-kirk:at}
James~F. Davis and Paul Kirk, \emph{Lecture notes in algebraic topology},
  Graduate Studies in Mathematics, vol.~35, American Mathematical Society,
  Providence, RI, 2001. \MR{\bf 1841974 (2002f:55001)}

\bibitem{finashin-kreck-viro:surfaces}
S.~M. Finashin, M.~Kreck, and O.~Ya. Viro, \emph{Nondiffeomorphic but
  homeomorphic knottings of surfaces in the {$4$}-sphere}, Topology and
  geometry---Rohlin Seminar, Lecture Notes in Math., vol. 1346, Springer,
  Berlin, 1988, pp.~157--198. \MR{\bf 970078 (90h:57021)}

\bibitem{finashin:surfaces}
Sergey Finashin, \emph{Knotting of algebraic curves in
  {$\mathbb{C}\mathrm{P}\sp 2$}}, Topology \textbf{41} (2002), no.~1, 47--55.
  \MR{\bf 1871240 (2003c:57042)}

\bibitem{fs:surfaces}
Ronald Fintushel and Ronald~J. Stern, \emph{Surfaces in {$4$}-manifolds}, Math.
  Res. Lett. \textbf{4} (1997), no.~6, 907--914. \MR{\bf 1492129 (98k:57047)}

\bibitem{fs:knots}
\bysame, \emph{Knots, links, and $4$-manifolds}, Invent. Math. \textbf{134}
  (1998), no.~2, 363--400. \MR{99j:57033}

\bibitem{fox:rolling}
R.~H. Fox, \emph{Rolling}, Bull. Amer. Math. Soc. \textbf{72} (1966), 162--164.
  \MR{\bf 0184221 (32 \#1694)}

\bibitem{freedman:concordance}
Michael~H. Freedman, \emph{A surgery sequence in dimension four;\ the relations
  with knot concordance}, Invent. Math. \textbf{68} (1982), no.~2, 195--226.
  \MR{\bf 666159 (84e:57006)}

\bibitem{freedman-quinn}
Michael~H. Freedman and Frank Quinn, \emph{Topology of 4-manifolds}, Princeton
  Mathematical Series, vol.~39, Princeton University Press, Princeton, NJ,
  1990. \MR{\bf 1201584 (94b:57021)}

\bibitem{GAP4}
The GAP~Group, \emph{{GAP -- Groups, Algorithms, and Programming, Version
  4.4}}, 2006, \verb+(http://www.gap-system.org)+.

\bibitem{hambleton-kreck:finite}
Ian Hambleton and Matthias Kreck, \emph{On the classification of topological
  {$4$}-manifolds with finite fundamental group}, Math. Ann. \textbf{280}
  (1988), no.~1, 85--104. \MR{\bf 928299 (89g:57020)}

\bibitem{hambleton-kreck:elliptic}
\bysame, \emph{Cancellation, elliptic surfaces and the topology of certain
  four-manifolds}, J. Reine Angew. Math. \textbf{444} (1993), 79--100.
  \MR{\bf 1241794 (95h:57036)}

\bibitem{hambleton-kreck:cancellation}
\bysame, \emph{Cancellation of hyperbolic forms and topological
  four-manifolds}, J. Reine Angew. Math. \textbf{443} (1993), 21--47.
  \MR{\bf 1241127 (94k:57030)}

\bibitem{hambleton-taylor:guide}
Ian Hambleton and Laurence~R. Taylor, \emph{A guide to the calculation of the
  surgery obstruction groups for finite groups}, Surveys on surgery theory,
  Vol. 1, Ann. of Math. Stud., vol. 145, Princeton Univ. Press, Princeton, NJ,
  2000, URL:http://www.math.mcmaster.ca/ian/published/Lgroups\_2000.pdf,
  pp.~225--274. \MR{\bf 1747537 (2001e:19007)}

\bibitem{kim:surfaces}
Hee~Jung Kim, \emph{Modifying surfaces in 4-manifolds by twist spinning}, Geom.
  Topol. \textbf{10} (2006), 27--56 (electronic). \MR{\bf 2207789}

\bibitem{kirby-siebenmann}
R.~C. Kirby and L.~C. Siebenmann, ``Foundational essays on topological
  manifolds, smoothings, and triangulations'', Princeton University Press,
  Princeton, N.J., 1977.
\newblock With notes by John Milnor and Michael Atiyah, Annals of Mathematics
  Studies, No. 88.

\bibitem{kirby-taylor:surgery}
Robion~C. Kirby and Laurence~R. Taylor, \emph{A survey of 4-manifolds through
  the eyes of surgery}, Surveys on surgery theory, Vol. 2 (S.~Cappell,
  A.~Ranicki, and J.~Rosenberg, eds.), Ann. of Math. Stud., vol. 149, Princeton
  Univ. Press, Princeton, NJ, 2001, pp.~387--421. \MR{2002a:57028}

\bibitem{kreck:surfaces}
Matthias Kreck, \emph{On the homeomorphism classification of smooth knotted
  surfaces in the {$4$}-sphere}, Geometry of low-dimensional manifolds, 1
  (Durham, 1989), London Math. Soc. Lecture Note Ser., vol. 150, Cambridge
  Univ. Press, Cambridge, 1990, pp.~63--72. \MR{\bf 1171891 (93k:57042)}

\bibitem{lawson:inertial}
Terry~C. Lawson, \emph{Inertial {$h$}-cobordisms with finite cyclic fundamental
  group}, Proc. Amer. Math. Soc. \textbf{44} (1974), 492--496. \MR{\bf 0358820
  (50 \#11279)}

\bibitem{litherland:deform}
R.~A. Litherland, \emph{Deforming twist-spun knots}, Trans. Amer. Math. Soc.
  \textbf{250} (1979), 311--331. \MR{\bf 530058 (80i:57015)}

\bibitem{madsen-milgram:surgery}
I.~Madsen and R.~J. Milgram, ``The classifying spaces for surgery and cobordism
  of manifolds'', vol.~92 of Annals of Mathematics Studies, Princeton
  University Press, Princeton, N.J., 1979.

\bibitem{montesinos:twins.II}
Jos{\'e}~M. Montesinos, \emph{On twins in the four-sphere. {I}}, Quart. J.
  Math. Oxford Ser. (2) \textbf{34} (1983), no.~134, 171--199. \MR{\bf 698205
  (86i:57025a)}

\bibitem{montesinos:twins.I}
Jos{\'e}~Mar{\'{\i}}a Montesinos, \emph{On twins in the four-sphere. {II}},
  Quart. J. Math. Oxford Ser. (2) \textbf{35} (1984), no.~137, 73--83.
  \MR{\bf 734666 (86i:57025b)}

\bibitem{morgan-bass:smith}
John~W. Morgan and Hyman Bass (eds.), \emph{The {S}mith conjecture}, Pure and
  Applied Mathematics, vol. 112, Academic Press Inc., Orlando, FL, 1984, Papers
  presented at the symposium held at Columbia University, New York, 1979.
  \MR{\bf 758459 (86i:57002)}

\bibitem{perron:isotopy2}
B.~Perron, \emph{Pseudo-isotopies et isotopies en dimension quatre dans la
  cat\'egorie topologique}, Topology \textbf{25} (1986), no.~4, 381--397.
  \MR{\bf 862426 (89g:57024)}

\bibitem{plotnick:fibered}
Steven~P. Plotnick, \emph{Fibered knots in {$S\sp 4$}---twisting, spinning,
  rolling, surgery, and branching}, Four-manifold theory (Durham, N.H., 1982),
  Contemp. Math., vol.~35, Amer. Math. Soc., Providence, RI, 1984,
  pp.~437--459. \MR{\bf 780592 (87a:57021)}

\bibitem{quinn:isotopy}
Frank Quinn, \emph{Isotopy of {$4$}-manifolds}, J. Differential Geom.
  \textbf{24} (1986), no.~3, 343--372. \MR{\bf 868975 (88f:57020)}

\bibitem{rourke-sullivan:arf}
C.~P. Rourke and D.~P. Sullivan, {\em On the {K}ervaire obstruction}, Ann. of
  Math. (2), {\bf 94} (1971), 397--413.

\bibitem{sullivan:gt}
Dennis~P. Sullivan, \emph{Geometric topology: localization, periodicity and
  {G}alois symmetry}, $K$-Monographs in Mathematics, vol.~8, Springer,
  Dordrecht, 2005, The 1970 MIT notes, Edited and with a preface by Andrew
  Ranicki. \MR{\bf 2162361}

\bibitem{sullivan:hauptvermutung}
D.~P. Sullivan, {\em Triangulating and smoothing homotopy equivalences and
  homeomorphisms. {G}eometric {T}opology {S}eminar {N}otes}, in ``The
  Hauptvermutung book'', vol.~1 of $K$-Monogr. Math., Kluwer Acad. Publ.,
  Dordrecht, 1996, 69--103.

\bibitem{wall:book}
C.~T.~C. Wall, \emph{Surgery on compact manifolds}, Academic Press, London,
  1970, London Mathematical Society Monographs, No. 1. \MR{\bf 0431216 (55
  \#4217)}

\bibitem{zeeman:twist}
E.~C. Zeeman, \emph{Twisting spun knots}, Trans. Amer. Math. Soc. \textbf{115}
  (1965), 471--495. \MR{\bf 0195085 (33 \#3290)}

\end{thebibliography}
%\bibliographystyle{amsplain}
\providecommand{\bysame}{\leavevmode\hbox to3em{\hrulefill}\thinspace}
\providecommand{\MR}{\relax\ifhmode\unskip\space\fi MR }
\providecommand{\MRhref}[2]{%
  \href{http://www.ams.org/mathscinet-getitem?mr=#1}{#2}
}
\providecommand{\href}[2]{#2}

\end{document}